\documentclass[11pt]{amsart}
\usepackage{anyfontsize}
\usepackage{tikz-cd}
\usepackage[headings]{fullpage}
\usepackage{amsmath}
\usepackage{amsfonts}
\usepackage{latexsym}
\usepackage{amssymb}
\usepackage{placeins}

\newtheorem{proposition}{Proposition}[section]
\newtheorem{theorem}{Theorem}[section]

\def\homo{\mathop{\sf Hom}}
\def\ker{\mathop{\sf Ker}}

\def\dim{\mathop{\sf dim}}

\def\sig{\sigma, \sigma^{-1}}

\def\B ellk{\mathcal{B}_{\ell,k}}
\def\S{\mathcal{S}}

\def\A{\mathcal{A}}
\def\C{\mathbb{C}}

\def\R{\mathbb{R}}
\def\Z{\mathbb{Z}}

\def\aut{\mathop{\sf Aut}}
\def\gl{\mathop{\sf GL_n(\mathbb{Z})}}
\def\glo{\mathop{\sf GL_1(\mathbb{Z})}}
\def\gla{\mathop{\sf GL_\mathbb{C}}}
\def\glt{\mathop{\sf GL_2(\mathbb{Z})}}

\newcommand{\pde}{[\partial_1, \ldots ,\partial_n]}

\newcommand{\diff}{[\sigma_1, \sigma_1^{-1}, \ldots ,\sigma_n, \sigma_n^{-1}]}

\newcommand{\F}{\mathcal{F}}

\begin{document}
 
\title[symmetric solutions] {Symmetric Solutions To Symmetric Partial Difference Equations}

\author{Shiva Shankar} 
\address{Department of Electrical Engineering, IIT Bombay, India.} \email{shunyashankar@gmail.com }

\begin{abstract} This paper studies discrete systems defined by linear difference equations on the lattice $\Z^n$ that are invariant under a finite group of symmetries, and shows that there always exist solutions to such systems that are also invariant under this group of symmetries.

\vspace{1mm}
\noindent {\bf Note}: Theorem 4.2 in [Symmetric Solutions to Symmetric Partial Difference Equations, ~Journal of Dfference Equations and Applications, ~31:3, 406-417, 2025, ~DOI:
10.1080/10236198.2024.2428380] does not charaterise  ~{\it all} \phantom{i}the symmetric solutions to a symmetric system of difference equations on $\Z^n,$ as claimed. Thus, the statement on the dimension of the space of symmetric solutions is an underestimate. 

\vspace{.7mm}
The correct statement appears in this arXiv version. Additionally, some typos have been corrected, and a few improvements incorporated.
 
\end{abstract}

\maketitle

\phantom{XXx} {\tiny Key words: linear partial difference equations, symmetry group, symmetric solution}

\phantom{XXx} {\tiny AMS classification: 39A05, 39A14, 58D19 }

\section{introduction}
This paper uses elementary 
methods from the representation theory of finite groups (acting on infinite dimensional vector spaces) 
to study solutions to a system of equations determined by a matrix of partial difference operators on the lattice $\Z^n$ that are invariant under a finite group of symmetries. It shows that if there exists a nonzero solution to the system of difference equations, then there also exists a nonzero symmetric solution, i.e. a solution that is invariant under this group of symmetries.

On the lattice $\Z^n$, denoted henceforth by $\mathbb{L}$, are $n$ independent forward shift operators $\sigma_1, \ldots, \sigma_n$ along the $n$ coordinate axes. These operators together with their inverses, the backward shifts, generate the algebra $A = \C\diff$ of partial difference operators on $\mathbb{L}$. Difference operators act on the (infinite dimensional) space $\F$ of all complex valued functions on $\mathbb{L}$. Given an $\ell \times k$ matrix $P(\sig)$ with entries from $A$, it defines a system of $\ell$ difference equations, and its solutions in $\F^k$ is the object of our study (precise definitions appear in the next section). 

Suppose that $G$ is a group of symmetries of the $\C$-algebra $A$; thus every element $g$ of $G$ is an algebra automorphism of $A$, and hence also a $\C$-vector space map. The group acts on $A^k$ component wise, and though this action does not define $A$-module maps, submodules of $A^k$ are still mapped to submodules by every $g$. Thus, given the $A$-submodule $P \subset A^k$ generated by the  $\ell$ rows of $P(\sig)$, it makes sense to ask if $P$ is $G$-invariant, i.e. if $g(P) = P$ for all $g$ in $G$. If so, then the system of equations defined by $P(\sig)$ is said to admit $G$ as a group of symmetries. Then, every $g$ in $G$ maps the set of solutions of $P(\sig)$ to itself (this is exactly similar to the case of a Lie point symmetry of a vector field). A symmetric solution to a $G$-invariant system is one which is invariant under the action of $G$.

Symmetry considerations in the study of differential equations is a subject with a long distinguished tradition, and is of singular importance in Physics. The invention of the spherical coordinate system to solve differential equations exhibiting rotational symmetry about a point (such as Laplace's equation), or of the cylindrical coordinate system to solve differential equations admitting rotational symmetry about an axis (as in Maxwell's equations describing the electromagnetic field produced by current flowing in a straight wire), are fundamental. E.~Noether's theorem in the setting of Lagrangian systems reduces the existence of first integrals to one prameter groups of symmetries of the Lagrangian \cite{a}. Thus the very possibility of solving a differential equation by quadrature is related to the existence of symmetries.

It is an important question to determine whether a system described by partial differential operators admits a vector potential \cite{pom}, a question also important in the setting of difference operators. This possibility is crucial
in the construction of wave-like trajectories of the system \cite{hhs}, as well as in engineering problems such as controller design \cite{stek}. An important issue here is the construction of these wave-like trajectories or controllers which are also invariant under the same group of symmetries as the system, and the answer reduces to the construction of a symmetric vector potential. Such questions will however not be pursued here.

While the subject of symmetries of differential equations is classical, the considerations of this article are specific to symmetries of difference equations. There is considerable work in this area, especially on the relation between symmetry and conservation laws, for instance \cite{lp}. Many papers deal with symmetries of a specific equation, in the setting of nonlinear difference equations (as in \cite{mfg}), as well as their boundary value problems. The departure here is that while this paper restricts itself to linear equations, it considers equations defined by matrix operators (i.e. submodules of $A^k$), and is different in nature from previous work in the subject. For instance, the paper of Fagnani and Willems \cite{fw} deals with symmetries of 1-$D$ systems on $\Z$, but its purpose is to obtain appropriate generators for an invariant system. Some of these results have been generalised in \cite{debs} to systems of partial difference equations on $\Z^n$ for a specific class of symmetries, but neither of these papers deals explicitely with symmetric solutions. The existence of such solutions naturally leads to finer questions on perturbations of discrete systems regarding the persistence of symmetric solutions when a symmetric system is perturbed within the class of symmetric systems. This would be the case especially when all the solutions to a symmetric system are also symmetric (Proposition 5.1 below).

The results here rest on the following correspondence, which is specific to difference equations. An element in the algebra $A$ is a sum of monomials; thus the set of monomials with coefficient 1 is a vector space basis for $A$. Such monomials also correspond to points of the lattice, the monomial $\sigma_1^{x_1} \cdots \sigma_n^{x_n}$ corresponding to the lattice point $(x_1, \ldots, x_n)$. This correspondence implies that the space $\F$ of complex valued functions on $\mathbb{L}$, in which solutions to difference equations are located, can be identified with the vector space dual $A'$ of $A$. This in turn allows G to act on $\F$, whence to the notion of a $G$-symmetric function, and finally to the notion of a symmetric solution to a symmetric system of difference equations.
\vspace{1.5mm}

Precisely then, the question addressed by this paper is the following: 
{\em if the system of difference equations is invariant under a finite group $G$ of symmetries of the ring $A$, then does the system admit solutions that are also $G$-invariant?}
\vspace{1mm}

The answer turns out to be in the affirmative. Furthermore, the results here provide a complete description of the space of $G$-invariant solutions, including an estimate of its dimension. This is possible only because the automorphism group of $A$ can be determined (Section 3). On the other hand, the automorphism group of the ring $\C\pde$ of partial differential operators on $\R^n$ is unknown \cite{d} (viz. the Jacobian Conjecture), and the details provided here for difference equations cannot be obtained in the setting of partial differential equations. \\

The following elementary example illustrates the nature of the results of this paper (details appear in the last section). 

\vspace{2mm}

Let $n = 1$ so that $\mathbb{L} = \Z$, and $A = \C[\sigma, \sigma^{-1}]$ is the ring of difference operators on $\mathbb{L}$ generated by the shift operator $\sigma$ and its inverse. Let the group $G = \Z/2\Z = \{1, -1\}$ act on the algebra $A$ by $-1(\sigma) = \sigma^{-1}$. The correspondence between monomials in $A$ and points in $\mathbb{L}$ alluded to above, identifies $\sigma^x$ with the point $x \in \mathbb{L}$, and hence a function $f: \mathbb{L} \rightarrow \C$ is invariant under the action of $G$ if $f(x) = f(-x)$, for all $x \in \mathbb{L}$. The space of such symmetric functions on $\mathbb{L}$ is thus infinite dimensional. (On the other hand, if the action on $A$ were to be given by $-1(\sigma) = -\sigma$, then an $f$ would be invariant if $f(x) = (-1)^xf(x)$, and hence if $f(2x+1) = 0$, for all $x \in \mathbb{L}$. This again defines an infinite dimensional space of functions on $\mathbb{L}$. There is yet another action, defined by $-1(\sigma) = -\sigma^{-1}$; it is considered in Example 3 of Section 5 below.)

(i) Consider the difference equation defined by the ideal $(\sigma + \sigma^{-1}) \subset A$; a solution $f: \mathbb{L} \rightarrow \C$ satisfies $\sigma(f)(x) + \sigma^{-1}(f)(x) = f(x+1) + f(x-1) = 0$. The ideal is invariant under the action of $G$, and a $G$-invariant solution is the sequence determined by $f(0) = 1, f(1) = 0$, namely the sequence $\cdots 0, -1, 0, \hat{1}, 0, -1, 0 \cdots$ (where $\hat{\phantom{x}}$ denotes the value of $f$ at $0 \in \mathbb{L}$). Indeed, the above equation admits a one dimensional space of $G$-invariant solutions, and the above $f$ spans this space.

(ii) Similarly, the difference equation defined by the ideal $(\sigma - \sigma^{-1})$ is also $G$-invariant. There is now a two dimensional space of $G$-invariant solutions. These are spanned by $\cdots 0, 1, 0, \hat{1}, 0, 1, 0 \cdots$ \phantom{x} and its translate \phantom{x} $\cdots 1, 0, 1, \hat{0}, 1, 0, 1, \cdots$. \\

The paper is organised as follows. The next section introduces notation and describes the structure of the set of solutions to a $G$-symmetric system of equations as a representation space for $G$. The third section describes the group $\aut_{\C}(A)$ of automorphisms of $A$. 
Section 4 proves the existence of symmetric solutions to systems of difference equations invariant under finite subgroups of $\aut_{\C}(A)$, and provides a concrete description of these solutions. The last section contains a few examples to illustrate these results, including the calculation of the dimension of the space of symmetric solutions. 

\section{solutions to a difference equation as a $G$-space} 

As in the introduction, let $\mathbb{L}$ denote the lattice $\Z^n = \{(x_1,\ldots, x_n) ~|~ x_i \in \Z, \forall i \}$ of all points in $\R^n$ with integral coordinates. For $i = 1, \ldots, n$, let $\sigma_i:\mathbb{L} \rightarrow \mathbb{L}$ map $x = (x_1, \ldots, x_i, \ldots, x_n)$ to $(x_1,\ldots, x_i+1, \ldots, x_n)$; it is the forward shift operator in the $i$-th direction. A monomial $\sigma^{x'} = \sigma_1^{x'_1}\cdots \sigma_n^{x'_n}$, $x'_i \in \Z$ for all $i$, maps points of $\mathbb{L}$ by composition; thus $\sigma^{x'}(x) = (x_1 + x'_1, \ldots, x_n + x'_n)$. This is the notation and terminology of \cite{debs}.

Let $A = \C\diff$ be the Laurent polynomial ring generated by these shifts, and their inverses. $A$ is a $\C$-algebra, it is the ring of partial difference operators on $\mathbb{L}$. Let $\F$ denote the set $\C^{\Z^n}$ of all complex valued functions on $\mathbb{L}$. The operator $\sigma_i$ induces an action on  $\F$ by  mapping $f \in \F$ to $\sigma_if$, where $\sigma_if(x) = f(\sigma_i(x))$. A monomial in $A$ acts by composition, and this action extends to an action of $A$ on $\F$, and  gives it the structure of an $A$-module.  

The object of study here is the set of solutions to a system of constant coefficient partial difference equations. Thus, if $P(\sig)$ is an $\ell \times k$ matrix whose entries $p_{ij}$ are in $A$, then it defines an $A$-module map
\[
\begin{array}{cccc}
P(\sig): & \F^k  & \longrightarrow & \F^\ell\\
& f=(f_1, \ldots, f_k) & \mapsto &  (p_1f, \ldots, p_\ell f) ,
\end{array}
\]
where the $i$-th row $p_i = (p_{i1}, \ldots, p_{ik})$ of $P(\sig)$ acts on $f$ by $p_if = \sum_{j=1}^kp_{ij}f_j$, $i=1, \ldots, \ell$. The set of solutions to this system of equations, namely the kernel $\ker_\F(P(\sig))$  of this map, depends on the $A$-submodule $P$ of $A^k$ generated by the rows of the matrix  $P(\sig)$, and not on the matrix itself. 
Indeed, the above kernel is isomorphic to $\homo_{A}(A^k/P, ~\F)$, the isomorphism given by the map
\[
\begin{array}{ccc}
\ker_\F(P(\sig))  & \longrightarrow & \homo_A(A^k/P, ~\F)  \\ 
f=(f_1, \ldots, f_k) & \mapsto & \phi_f   ~,
\end{array}
\]
where $\phi_f([(q_1, \ldots, q_k)]) = q_1f_1 + \cdots + q_kf_k$ (the square brackets denotes the class of an element of $A^k$ in the quotient $A^k/P$). Hence this kernel will be denoted by $\ker_\F(P)$; it is the set of solutions in $\F^k$ to the system of difference equations defined by the submodule $P$ of $A^k$.

The results of this paper stem from the following observation. If the point $x = (x_1, \ldots ,x_n)$ in $\mathbb{L}$ is identified with the monomial $\sigma^x = \sigma_1^{x_1} \cdots \sigma_n^{x_n}$ in $A$, then $A$ can be identified with the $\C$-vector space spanned independently by the points of $\mathbb{L}$. The space $\F$ is then isomorphic to the vector space dual $A' = \homo_\C(A, ~\C)$ of $A$. Given elements $f \in A'$ and $a \in A$, define $af \in A'$ by $af(a') = f(aa')$.\footnote{Henceforth $f$ will denote either a function $f: \mathbb{L} \rightarrow \C$ or its extension to a linear map $f:A \rightarrow \C$.} This gives $A'$  the structure of an $A$-module, and the above $\C$-isomorphism is an isomorphism of $A$-modules.
Similarly, $\homo_\C(A^k/P, ~\C)$ is an $A$-module, and it follows that\[{\ker}_\F(P) \simeq {\homo}_A(A^k/P, ~{\homo}_\C(A, ~\C))  \simeq {\homo}_\C(A^k/P \otimes_A A, ~\C) \simeq {\homo}_\C(A^k/P, ~\C)\]by tensor-hom adjunction. Thus the set of solutions to the system $P$ is isomorphic to the dual $(A^k/P)'$ as $A$-modules, where an $f = (f_1, \ldots, f_k) \in \ker_\F(P)$ is identified with $f: A^k/P \rightarrow \C$ mapping $[(q_1, \ldots, q_k)]$ to $f_1(q_1) + \cdots + f_k(q_k)$, each $f_i \in \F$ now defined on $A$ by linear extension.

It is this identification that allows the use of group representation methods in the study of solutions to partial difference equations. \\

\noindent Remark 2.1. ~As the functor $\homo_\C(-, ~\C)$ is exact on the category of $\C$-vector spaces, so is the functor $\homo_A(-,~\F)$ on the category of $A$-modules. This implies that $\F$ is an injective $A$-module. Moreover, $\homo_A(M, ~\F) \neq 0$, if $M \neq 0$, and thus $\F$ is an injective cogenerator.

The injective cogenerator property of $\F$ implies that the assignment $P \rightarrow \ker_\F(P)$ is an inclusion reversing  bijective correspondence between $A$-submodules of $A^k$ and its solutions in $\F^k$, see for instance \cite{stek}. \\  

Let $(A, \rho)$, $\rho: G \rightarrow \gla(A)$, be a representation of a group $G$ over $\C$. In turn, $G$ acts on $\F$ by $g \cdot f(x) = f(\rho(g^{-1})x)$ ~- ~this is the dual representation $\rho'$ of $G$ on $A'$ ~- ~and to say that the function $f$ is symmetric or invariant with respect to $G$ is to say that $f$ is a fixed point for $\rho'$ and hence constant on an orbit of $G$ in $A$. If $G$ is a subgroup of the group $\aut_\C(A)$ of $\C$-algebra automorphisms of $A$, then $G$ maps monomials bijectively to monomials (Section 3 below). The correspondence between monomials with coefficient 1 and points of the lattice $\mathbb{L}$ then implies that if $\rho(g)x = \alpha y$, $x, y \in \mathbb{L}, ~\alpha \in \C^*$, then a $G$-invariant $f$ satisfies $f(x) = \alpha f(y)$. In particular, if $G$ maps the set of monomials with coefficient 1 bijectively to itself, so that $\rho$ restricts to an action of $G$ on the set $\mathbb{L}$, then a $G$-invariant $f$ assumes the same value on the points of $\mathbb{L}$ that are in the same $G$-orbit.

Given the above representations $\rho$ and $\rho'$ of $G \subset \aut_\C(A)$, they define the direct sum representations of $G$ on $A^k$ and on $(A^k)'$, again denoted $\rho$ and $\rho'$ respectively. 
If a system of difference equations given by a submodule $P \subset A^k$ is left invariant by $G$, i.e. if $\rho(g)P = P$ for all $g \in G$, then $G$ acts on $A^k/P$, and hence on $(A^k/P)' \simeq \ker_\F(P)$. A fixed point $f$ of this representation is a symmetric or invariant solution of $P$, it satisfies $f(x) = 
f(\rho(g^{-1})x)$, for all $x \in \mathbb{L}$ and for all $g \in G$. (Henceforth, we suppress $\rho$ and write $gx$ for $\rho(g)x$.)

The rest of the paper assumes that $G$ is a {\it finite} subgroup of $\aut_\C(A)$. This allows the use of results from the theory of finite dimensional representations of $G$ even though $A^k/P$ might not be finite dimensional. \\

The following facts about finite dimensional representations of a finite group $G$ carry over to infinite dimensional representations. 
\vspace{1mm}

Let $V$ be a representation of $G$. Averaging over $G$ is possible, and hence if $W_1 \subset V$ is a subrepresentation, then there is a subrepresentation $W_2 \subset V$ such that $V = W_1 \bigoplus W_2$. Irreducible representations of $G$ are finite dimensional, and there are finitely many of them, say $V_1, \ldots, V_t$ (equal to the number of conjugacy classes of $G$). Representations are completely reducible, and the proof is by induction as in the case of a finite dimensional representation. 
Thus, let $\frak{P} = \{V = \bigoplus\limits_{R_1}V_1 \cdots \bigoplus\limits_{R_t}V_t \bigoplus W ~|~ W ~is ~ G$-$invariant\}$.
Partially order $\frak{P}$ by setting $\bigoplus\limits_{R_1}V_1 \cdots \bigoplus\limits_{R_t}V_t \bigoplus W \leqslant \bigoplus\limits_{R'_1}V_1 \cdots \bigoplus\limits_{R'_t}V_t \bigoplus U$ if $U \subseteq W$. Every chain in $\frak{P}$ has an upper bound, and hence there exists a maximal element in $\frak{P}$ which must be of the form $\bigoplus\limits_{S_1}V_1 \cdots \bigoplus\limits_{S_t}V_t$. This is the required decomposition of $V$ as a direct sum of irreducibles. 
Furthermore, the decomposition is unique.

It now follows that 
\begin{equation}
V' = \prod\limits_{S_1} \homo{_\C}(V_1, \C) \phantom{.}\cdots \phantom{.} \prod\limits_{S_t} \homo{_\C}(V_t, \C),
\end{equation}
which is written briefly as $\prod\limits_{i,S_i}V_i'$. As each $V_i$ is an irreducible finite dimensional representation of $G$, so is its dual, and (1) expresses $V'$ as a product of irreducible representations.

In particular, the solutions to a system of $G$-invariant partial difference equations can be expressed as a product of irreducible representations.

\begin{proposition} Let $G$ be a finite subgroup of $\aut_\C(A)$ and let $V_1, \ldots, V_t$ be its irreducible representations ($V_1$ the trivial representation). Let $P$ be a submodule of $A^k$ which is invariant for the direct sum representation of $G$ on $A^k$. 
Then $A^k/P$ can be expressed uniquely as $\bigoplus\limits_{S_1}V_1 \cdots \bigoplus\limits_{S_t}V_t$.
Hence, the set of solutions in $\F$ of the system of difference equations defined by $P$ admits the product decomposition
\[
\ker{_\F}(P) \simeq (A^k/P)' \simeq \prod_{S_1}V_1' \phantom{.} \cdots \phantom{.} \prod_{S_t}V_t'.
\]		
\end{proposition}

Symmetric solutions to the system of difference equations defined by $P$ are the fixed points $\ker_\F(P)^G$ of the $G$-action on $\ker_\F(P)$. The article turns to the description of these fixed points in Section 4 after describing the finite subgroups of $\aut_\C(A)$.

\section{finite subgroups of $\aut_\C(A)$}

Recollect first the calculation of the automorphism group $\aut_{\C}(A)$ of $A$ from \cite{debs}.
A $\C$-algebra endomorphism of $A$ maps $\sigma_i$ to a unit, and hence to a monomial $r_i\sigma_1^{m_{i1}} \cdots \sigma_n^{m_{in}}$, where the $r_i \in \C^*$ and the $m_{ij} \in \Z$, for $i, j = 1, \ldots, n$.
If this endomorphism is a $\C$-algebra automorphism, then the $n \times n$ matrix $(m_{ij})$ is unimodular.
Thus, a general automorphism of $A$ is obtained from two 
$\C$-algebra automorphisms, the first given by a homothety, $\sigma_i \mapsto r_i\sigma_i$, and the second given by $\sigma_i \mapsto \sigma_1^{m_{i1}} \cdots \sigma_n^{m_{in}}$, for all $i$. They define two group homomorphisms:
\vspace{1mm}

(i) $\psi_1: (\C^*)^n \rightarrow \aut_\C(A)$, defined by $\psi_1(R)\sigma_i = r_i\sigma_i$, $i = 1, \ldots, n$, where $R = (r_1, \ldots, r_n)$,

(ii) $\psi_2: \gl \rightarrow \aut_\C(A)$, defined by $\psi_2(M)\sigma_i = \sigma_1^{m_{i1}} \cdots \sigma_n^{m_{in}}$, $i = 1, \ldots, n$, where the entries of $M$ are the $m_{ij}$. 
\vspace{1mm}

Let $\phi: \gl \rightarrow  \aut((\C^*)^n)  $ be the group homomorphism defined by $\phi(M)(R) = (\prod_{i =1}^n r_i^{m_{1i}}, \ldots, \prod_{i=1}^n r_i^{m_{ni}})$, where $R = (r_1,\ldots, r_n)$ and $M$ is the matrix $(m_{ij})$. Let $(\C^*)^n \rtimes_\phi \gl$ be the semi-direct product of $(\C^*)^n$ and $\gl$ determined by $\phi$. Then
\[ \psi_1(R)\circ \psi_2(M)\sigma_i = \psi_2(M) \circ  \psi_1(\phi(M)(R))\sigma_i, \]
for all $i$, and for all $R \in (\C^*)^n, M \in \gl$. Thus, the actions defined by $\psi_1$ and $\psi_2$ lift to the semidirect product, and gives a homomorphism $\Psi: (\C^*)^n \rtimes_\phi \gl \rightarrow \aut_\C(A)$.

Conversely, as an element $g$ of $\aut_\C(A)$ is determined by the images of the $\sigma_i$ under $g$, it follows that if $g\sigma_i$ is equal to $r_i\sigma_1^{m_{i1}} \cdots \sigma_n^{m_{in}}$, $i = 1, \ldots, n$, then the $g\sigma_i$ determine the element $R = (r_1, \ldots, r_n) \in (\C^*)^n$ and the unimodular matrix $M = (m_{ij})$. This defines $\aut_\C(A) \rightarrow (\C^*)^n \rtimes_\phi \gl$, by mapping $g$ to $(R, M)$, which is inverse to $\Psi$.

This establishes the following proposition.

\begin{proposition} $\aut_\C(A) \simeq (\C^*)^n \rtimes_\phi \gl$. 
\end{proposition}

Both $(\C^*)^n$ and $\gl$ will be considered to be subgroups of $\aut_\C(A)$, the first via the canonical inclusion and the second via the splitting $\gl \stackrel{i'}{\hookrightarrow} (\C^*)^n \rtimes_\phi \gl$ of the exact sequence
$1 \rightarrow (\C^*)^n \stackrel{i}{\longrightarrow} (\C^*)^n \rtimes_\phi \gl \stackrel{p}{\longrightarrow} \gl \rightarrow 1.$ \\

A finite subgroup of $(\C^*)^n$ is a product of at most $n$ cyclic groups, and there are infinitely many such. On the other hand, there are only a finite number of finite subgroups of $\gl$, by a theorem of Minkowski. Subgroups of a semidirect product do not admit a characterization in general (such as Goursat's theorem for a direct product), but the results of the next section are not dependent on any specific description of finite subgroups of $\aut_\C(A)$.

\section{symmetric solutions}
This section establishes the existence of symmetric solutions to systems of difference equations admitting a finite group of symmetries, and describes its structure.

\begin{proposition} Let $G$ be a finite subgroup of $\aut_\C(A)$ and let $V_1, \ldots, V_t$ be its irreducible representations, where $V_1$ is the trivial representation. Let $P \subset A^k$ be invariant under $G$. Let $A^k/P = \bigoplus\limits_{S_1}V_1 \cdots \bigoplus\limits_{S_t}V_t$ be the decomposition into a sum of irreducibles. Then the dimension of the space $\ker_\F(P)^G$ of symmetric solutions to $P$ equals the cardinality of $S_1$ if finite, and is infinite dimensional otherwise.
\end{proposition}
\noindent Proof: As $V_1$ is the trivial representation, so is its dual $V_1'$, and hence $\dim(\homo(V_1, V'_i)^G)$ is 1 when $i = 1$ and is 0 otherwise. The dimension of $\ker_\F(P)^G$ equals the multiplicity of $V_1$ in $\ker_\F(P) = \prod\limits_{i,S_i}V_i'$ (in the notation of Section 2), and thus equals the dimension of $\homo(V_1, \prod\limits_{S_1}V_1')$. This proves the proposition.
\hspace*{\fill}$\square$\\

\begin{theorem} Let $P$ be a proper submodule of $A^k$ which is invariant under the action of a finite subgroup $G$ of $\aut_\C(A)$. Then there are nontrivial symmetric solutions to the difference equation defined by $P$.
\end{theorem}
\noindent Proof: Without loss of generality, it can be assumed that $e_1 = (1, 0, \ldots, 0) \notin P$. Then $G$ acts trivially on the 1-dimensional space spanned by the coset $[e_1]$ in $A^k/P$, and hence the multiplicity of the trivial representation $V_1$ in $A^k/P$ is at least 1. The theorem follows from the above proposition.
\hspace*{\fill}$\square$\\ 

{\em Thus, if $P$ admits nonzero solutions, then it also admits nonzero symmetric solutions.}

\vspace{2mm}

The rest of this section describes the space of symmetric solutions to a symmetric system guaranteed by the above theorem. It first considers finite subgroups of each of the two subgroups $(\C^*)^n$ and $\gl$ of $\aut_\C(A)$ separately, and then  describes symmetric solutions invariant under an arbitrary finite subgroup of $\aut_\C(A)$. (The description of symmetric solutions in the case of finite subgroups of $(\C^*)^n$ relies on results in \cite{debs}, but the results below are otherwise independent of it.) \\

Suppose first that $G$ is a finite subgroup of $(\C^*)^n \subset \aut_\C(A)$. 
An element $\zeta = (\zeta_1, \ldots, \zeta_n) \in G$, each $\zeta_i$ a root of unity, acts on a monomial $\sigma^x = \sigma_1^{x_1} \cdots \sigma_n^{x_n}$ to give $\zeta \sigma^x = (\zeta_1\sigma_1)^{x_1} \cdots (\zeta_n\sigma_n)^{x_n} = \zeta_1^{x_1} \cdots \zeta_n^{x_n}\sigma^x$ (namely Section 3), and thus an element of $A$ is invariant under the action of $G$ only if each of its monomial terms is. 
The space $A^G$ of all the elements of $A$ that are left fixed by every element of $G$ is a subalgebra of $A$. As $g(\sigma^x)g(\sigma^{-x}) = \sigma(1) = 1$, it follows that if $\sigma^x \in A^G$, then so is $\sigma^{-x}$. Thus, if a monomial is in $A^G$, its inverse is also in it, and $A^G$ is a Laurent subalgebra of $A$. The ring $A$ is a free $A^G$-module of finite rank equal to the cardinality $|G|$ of $G$.

It follows that the set of monomials in $A^G$ with leading coefficient 1 is a $\C$-vector space basis for it. Let $\mathbb{S}$ be the set of points of $\mathbb{L}$ corresponding to these monomials. It follows, as $A^G$ is a Laurent algebra, that $\mathbb{S}$ is a sublattice of $\mathbb{L}$, of full rank as $G$ is finite, whose ring of difference operators $A_S$ is precisely $A^G$. Let $\F_S$ denote the set of all complex valued functions on $\mathbb{S}$, it is an $A_S$ module. It is now clear that an $f: \mathbb{L} \rightarrow \C$ is $G$-invariant if and only if the support of $f$ is contained in $\mathbb{S}$, i.e. if $f(x) = 0$ for all $x \notin \mathbb{S}$. The subspace of $G$-invariant functions on $\mathbb{L}$ is isomorphic to $\F_S$ and is thus infinite dimensional. 

Let $P^c$ denote the contraction $P \cap A_S^k$ of $P$ to $A_S^k$.

\begin{proposition} Let $G$ be a finite subgroup of $(\C^*)^n \subset \aut_\C(A)$. Let $\mathbb{S}$ be the sublattice of ~$\mathbb{L}$ on which $G$ acts trivially, and $\F_S$ the space of functions on it. Let $P$ be a $G$-invariant submodule of $A^k$. Then $\ker_\F(P)^G \simeq \ker_{\F_S}(P^c)$, so that the space of symmetric solutions to $P$ is of the same dimension as the space of solutions to the contracted system $P^c$ on $\mathbb{S}$. This dimension is infinite unless $\dim(A^k/P)$ is finite, in which case $\dim(\ker_\F(P)^G) = \frac{1}{|G|} \dim(A^k/P)$.
\end{proposition}
\noindent Proof: The proof follows from an explicit description of $\ker_\F(P)$ in \cite{debs} for a $G$-invariant submodule $P \subset A^k$. Given $P^c$, the contraction of $P$ to $A^k_S$, let $P^{ce}$ denote the extension of $P^c$ back to $A^k$. To say that $P$ is $G$-invariant is equivalent to $P^{ce} = P$ (Proposition 5.3 of \cite{debs}). Let $\ker_{\F_S}(P^c) \simeq \homo_{A_S} (A_S^k/P^c, ~\F_S)$ be the solutions to the contracted system $P^c$ on $\mathbb{S}$. Then by Propositions 4.4 and 3.3 of \cite{debs}, $\ker_\F(P) \simeq \prod_{|G|}\ker_{\F_S}(P^c)$, where a solution to the system $P$ is constructed by choosing $|G|$ solutions to the contracted system $P^c$ on the sublattice $\mathbb{S}$, and assigning them to $\mathbb{S}$ and its translates (cosets) in $\mathbb{L}$.

Given this description, a $G$-invariant solution to $P$ is constructed by assigning an arbitrary solution of $P^c$ to $\mathbb{S}$ and the 0 solution to its translates in $\mathbb{L}$. Thus $\ker_\F(P)^G$ is isomorphic to $\ker_{\F_S}(P^c)$ as $\C$-vector spaces. This proves the proposition.
\hspace*{\fill}$\square$\\

Next suppose that $G$ is a finite subgroup of $\gl \subset \aut_\C(A)$. $G$ maps the set of monomials in $A$ with coefficient 1 bijectively to itself (Section 3), hence the action of $G$ on $A$ restricts to $\mathbb{L}$. This action is the `point symmetries' analogue of the introduction. 

Let $x \in \mathbb{L}$ and let $\mathbb{L}_x = ~<x, y, \ldots, z>$ be its $G$-orbit in $\mathbb{L}$ (the brackets $< ~>$ indicate that a point in the orbit might occur several times). Let $A_x$ be the subspace of $A$ spanned by $\mathbb{L}_x$, and let $(A_x, \rho_x)$ be the representation of $G$ defined by its action on $\mathbb{L}_x$. Let $\rho_x = \rho_1 \oplus \rho_2$ is its decomposition into a direct sum, where $\rho_1$ is the trivial representation on the subspace $A_{x'}$ of $A_x$ spanned by $x + y + \cdots + z$, and $\rho_2$ is the standard representation on the subspace $A_{x''}$ of $A_x$ given by $\{c_xx + \cdots + c_zz ~|~ c_x + \cdots + c_z = 0 \}$. Thus, $A_x^G$ is 1 dimensional, and so is the space of $G$-invariant functions in $A'_x$. It is spanned by any function which is nonzero  on $x + y + \cdots + z$ and 0 on $A_{x''}$, for instance $f$ such that $f(x)= f(y)= \cdots = f(z) = 1$

Hence if $\mathbb{L} = \bigsqcup_{i \in I} \mathbb{L}_i$ is the disjoint union of orbits, then the subrepresentation of $G$ defined on the subspace $A_i$ of $A$ spanned by the elements of $\mathbb{L}_i$ has a 1 dimensional space $A_i^G$ of fixed points, for each $i \in I$. As $A = \bigoplus_{i \in I}A_i$, it follows that $A^G = \bigoplus_{i \in I}A_i^G$ has dimension equal to the cardinality of $I$, and hence that there is an infinite dimensional subspace of $G$-invariant functions in $\F$. 

Similarly, given the direct sum representation of $G$ on $A^k$, the orbit $\mathbb{L}_i$ defines a subrepresentation $A_i^k$ and a $k$-dimensional space $(A_i^k)^G \simeq (A_i^G)^k$ of fixed points, for each $i$. 

Now let $P \subset A^k$ be $G$-invariant; then $P^G = P \cap (A^G)^k = P \cap \bigoplus_{i \in I} (A_i^G)^k$.

\begin{proposition} Let $G$ be a finite subgroup of $\gl \subset \aut_\C(A)$. Let $\mathbb{L} = \bigsqcup_{i \in I}\mathbb{L}_i$ be the disjoint union of orbits of the action of $G$, and $A^k = \bigoplus_{i \in I}A_i^k$ the corresponding direct sum of $G$-representations. Let $P$ be a $G$-invariant submodule of $A^k$. Then the space of symmetric solutions to $P$ is infinite dimensional unless there are only finitely many indices, say $I_f \subset I$, such that the projection $\pi: A^k \rightarrow \bigoplus_{i \notin I_f} A_i^k$ satisfies $\bigoplus_{i \notin I_f}(A_i^G)^k \subset \pi(P)$.	\\
\end{proposition}
\noindent Proof: By Propositions 4.1 and 2.1, there is no contribution to $\ker_\F(P)^G$ from the subspace of $(A_i^G)^k$ given by $\pi_i(P) \cap (A_i^G)^k$, where $\pi_i:A^k \rightarrow A_i^k$ is the projection to the $i$-th summand. Thus the contribution from $(A_i^G)^k$ is of dimension 0 if and only if $(A_i^G)^k \subset \pi_i(P)$. The proposition now follows.
\hspace*{\fill}$\square$\\


Consider finally the case of an arbitrary finite subgroup $G \subset \aut_\C(A)$.
For $x\in \mathbb{L}$, let $\mathbb{A}_x$ now denote its orbit in $A$ (in Proposition 4.3,  this orbit was contained in $\mathbb{L}$ and was denoted $\mathbb{L}_x$). Let $\A = \{\mathbb{A}_x ~|~ x \in \mathbb{L}\}$ be the collection of all the orbits through all the points in $\mathbb{L}$. An orbit $\mathbb{A}_x$ is {\it degenerate} if the sum of all the elements in it is equal to 0, and is {\it nondegenerate} otherwise. Let $\A_d$ and $\A_n$ denote the collection of degenerate and nondegenerate orbits respectively.

\begin{theorem}
 Let $G$ be a finite subgroup of $\aut_\C(A)$. Let $A = \bigoplus_{i \in \A_d} A_i ~\bigoplus_{i \in \A_n}A_i$ be the direct sum representation of $G$ corresponding to the disjoint union of orbits through points in $\mathbb{L}$. Then, $(A_i)^G$ is 1-dimensional for $i \in \A_n$ and 0-dimensional for $i \in \A_d$. Thus, each $i \in \A_n$ contributes $k$ to the dimension of $(A^k)^G \simeq \bigoplus_{i \in \A_n} (A_i^k)^G$. 
 
 Let $P$ be a $G$-invariant submodule of $A^k$. Then the space $\ker_\F(P)^G$ of all symmetric solutions to $P$ is infinite dimensional unless $\dim(A^k/P)$ is finite or there are finitely many $I_f \subset \A_n$ such that the projection $\pi: A^k \rightarrow \bigoplus_{i \notin I_f} A^k$ satisfies $\bigoplus_{i \notin I_f}(A_i^k)^G \subset \pi(P)$.
\end{theorem}
\noindent Proof: Let the orbit $\mathbb{A}_x$ through $x \in \mathbb{L}$ be denoted by $<x, \eta_yy, \ldots, \eta_zz>$, where $y, \ldots,z$ are monomials with leading coefficient 1, hence correspond to points in $\mathbb{L}$, and the coefficients $\eta_y,\ldots, \eta_z$ are in $\mathbb{C}^*$ (in Proposition 4.3, these coefficients are all 1 as $\mathbb{A}_x \subset \mathbb{L}$). 

Suppose first that $<x, \eta_yy, \ldots, \eta_zz>$ is linearly independent, and thus $\mathbb{A}_x \in \A_n$; this is exactly when  $x, y,\ldots, z$ are all distinct. Let $A_1$ be the subspace of $A_x$ spanned by $x + \eta_yy + \cdots +\eta_zz$, and $A_2$ the subspace given by $\{c_xx+c_y\eta_yy+ \cdots +c_z\eta_zz ~|~ c_x+c_y+ \cdots +c_z = 0\}$.
Then, exactly as in Proposition 4.3, the representation $(A_x,\rho_x)$ decomposes as a direct sum $(A_{x'}, \rho_1) \oplus (A_{x''}, \rho_2)$, where the first is the trivial representation, and the second is the standard representation on $A_{x''}$. Thus, $A_x^G$ is 1 dimensional, and so is the space of $G$-invariant functions in $A'_x$. It is spanned by any function which is nonzero  on $x + \eta_yy + \cdots +\eta_zz$, and 0 on $A_{x''}$, for instance $f$ such that $f(x)=1,~f(y)=\eta_y^{-1},\ldots, f(z) = \eta_z^{-1}$.  

\vspace{1.5mm}
Suppose next that $<x, \eta_yy, \ldots, \eta_zz>$ is linearly dependent, and hence that $x, y, \ldots, z$ are {\it not} distinct. Suppose that $x, \zeta_2x, \ldots, \zeta_rx$ are all the points in the orbit $\mathbb{A}_x$ that are dependent on $x$, where say $x = ex, ~\zeta_2x = g_2x, \ldots, \zeta_rx = g_rx$. It follows that $G(x) = \{e, g_2, \ldots, g_r\}$ is a subgroup of $G$, and that if $G_x$ is the stabilizer of $x$, then the quotient $G(x)/G_x$ is isomorphic to the group of $k$-th roots of unity, for some $k$. It follows that each $\zeta_ix$ occurs the same number of times in the $G(x)$-orbit $<x,\zeta_2x, \ldots, \zeta_rx>$. 

For any other point in the $G$-orbit of $x$, say $\eta_yy$, the subgroup $G(\eta_yy)$ of elements of $G$ that map $\eta_yy$ to a multiple of it, is conjugate to $G(x)$. Its stabilizer $G_{\eta_yy}$ is isomorphic to $G_x$, and hence it again follows that $G(\eta_yy)/G_{\eta_yy}$ is isomorphic to the same $k$-th roots of unity as above. Each element in the $G(\eta_yy)$-orbit of $\eta_yy$  thus occurs the same number of times as $\eta_yy$. 

The proof now splits into two cases:

\noindent (i) Suppose $G_x \simeq G(x)$, i.e. the point $x$ is $G(x)$-symmetric. Let $<x, \
\eta_uu,\ldots, \eta_vv>$ be a maximally independent set of points in the orbit $\mathbb{A}_x$. Then the sum of all the elements in $\mathbb{A}_x$ equals $k(x + \eta_uu + \cdots + \eta_vv)$, hence this orbit belongs to $\A_n$. An easy calculation shows that there is a unique 1-dimensional subspace of $A_x$, spanned by $x + \eta_uu + \cdots + \eta_vv$ which is $G$-invariant. The 1 dimensional space of $G$ invariant functions in $A'_x$ is constructed just as above: it is spanned by $f(x) = 1, ~f(u) = \eta_u^{-1}, \ldots,~f(v) = \eta_v^{-1}$. 

\noindent (ii) Suppose that $G_x\subsetneq G(x)$. Then the sum of the elements in the $G(x)$-orbit of $x$ equals a multiple of $\sum_{i=1}^k \zeta^i x$, where $\zeta$ is a $k$-th root of unity, and thus equal to 0. As this is true for every point in $\mathbb{A}_x$, the sum of all the points in it equals 0, hence $\mathbb{A}_x \in \A_d$. A routine calculation shows that there is no nontrivial $G$-invariant vector in $A_x$, and hence no nontrivial $G$-invariant function in $A'_x$. 

This proves the first part of the theorem. The second part follows just as in Proposition 4.3
\hspace*{\fill}$\square$\\

If $V$ is a finite dimensional complex representation of a finite group $G$, then the dimension of the subspace $V^G = \{v \in V ~| ~gv = v, \forall g \in G\}$ of fixed points is given by $\frac{1}{|G|}\sum_{g \in G} \chi_{\scriptscriptstyle V}(g)$, where $\chi_{\scriptscriptstyle V}$ is the character of $V$. However, here the representation $(A^k/P)'$ could be infinite dimensional and the above formula meaningless.

Instead, the discussion preceeding the statement of Proposition 4.3 can be adapted to provide precise information on the dimension of symmetric solutions. Thus, suppose now that $V = A^k/P$ is infinite dimensional. Exhaust $V$ by an increasing sequence of finite dimensional $G$-invariant subspaces $W_0 \subset \cdots W_i \subset W_{i+1} \subset \cdots $~. For each $i$, $\dim(W_i^G) = \frac{1}{|G|}\sum_{g \in G} \chi_{\scriptscriptstyle W_i}(g)$.  As $V^G = \bigcup_iW_i^G$, it follows that the dimension of $V^G$ is either infinite, when the sequence $\{\dim(W_i^G)\}$ is unbounded, or is equal to its limit, otherwise.  In the latter case, there is an index $i_0$ such that $W_i^G = W_{i_0}^G$ for all $i \geqslant i_0$. The value of $\dim(V^G)$ is clearly independent of the choice of the exhaustion $\{W_i\}$.

This proceedure is employed in the next section to explain the assertions in the examples of the introduction.

\section{examples}
This section consists of several examples that illustrate the above results. \\

\noindent 1. In the examples of the introduction, $A = \C[\sigma, \sigma^{-1}]$ is the ring of difference operators on the lattice $\mathbb{L} = \Z$, and $G = \{1, -1\}$ acts as $\C$-algebra automorphisms on $A$ by 
$-1(\sigma) = \sigma^{-1}$ ($G$ is the subgroup $\glo \subset \aut_\C(A)$). This is the group of reflections of $\Z$ about the origin.

Exhaust $A$ by the nested sequence of $G$-invariant finite dimensional subspaces $\{W_i\}_{i \geqslant 0}$, where $W_i$ is the subspace of dimension $2i + 1$ spanned by $\{\sigma^x,~ -i \leq x \leq i\}$. The map $-1: W_i \rightarrow W_i$ has trace 1, and hence the subspace $W_i^G$ of fixed points of $W_i$ has dimension $\frac{1}{2}(\chi(1) + \chi(-1)) = i + 1$. Thus, $A^G$ is infinite dimensional, and hence $A'^G$, the space of all $G$-invariant functions in $\F$, is also infinite dimensional.
\vspace{1.5mm}

(i) Let $I = (\sigma + \sigma^{-1})$, it is a $G$-invariant ideal of $A$. Then for $i \geqslant 1$, $I_i = I \cap W_i$ is a subspace of $W_i$ of dimension $2i-1$ (spanned by $\{\sigma^i+\sigma^{i-2}, \ldots, \sigma^{-i+2}+\sigma^{-i} \}$), hence the trace of $-1: I_i \rightarrow I_i$ equals 1. Thus $I_i^G$ has dimension $\frac{1}{2}(\chi(1) + \chi(-1)) = i$, and hence $\dim(W_i/I_i)^G = \dim(W_i^G) - \dim(I_i^G) = 1$, for all $i$. It follows that the space $(A/I)'^G$ of $G$-invariant solutions to the difference equation defined by $I$ has dimension 1 (and is spanned by the solution $f$ given by $\cdots 0, -1, 0, \hat{1}, 0, -1, 0 \cdots$, where, as before, $\hat{\phantom{x}}$ denotes the value of $f$ at 0.)

Indeed, in the notation of Proposition 4.3, $I_f$ contains a single element corresponding to the orbit $\{ 1 \}$, as is readily verified.

The dimension of $W_i/I_i$ equals 2 for all $i$, hence the dimension of $(A/I)'$ equals 2. The solutions of $I$ other than the span of $f$ above is spanned by its (non-symmetric) shift $\cdots 1, 0, -1, \hat{0}, 1, 0, -1, \cdots$. 
\vspace{1.5mm}

(ii) Let $J$ be the $G$-invariant ideal $(\sigma - \sigma^{-1})$; the dimension of the $G$-invariant subspace $J_i = J \cap W_i$ is also $2i-1$ for $i \geqslant 1$. The trace of $-1: J_i \rightarrow J_i$ is now $-1$, hence $\dim(J_i^G) = \frac{1}{2}(\chi(1) + \chi(-1)) = i - 1$. Then $\dim((W_i/J_i)^G) = 2$ for all $i$, and thus the dimension of the space $(A/J)'^G$ of $G$-invariant solutions to the equation defined by $J$ equals 2 (spanned by $\cdots 0, 1, 0, \hat{1}, 0, 1, 0 \cdots$ and $\cdots 1, 0, 1, \hat{0}, 1, 0, 1, \cdots$). 

The dimension of $(A/J)'$ equal 2, and hence {\it all} the solutions to the equation defined by $J$ are symmetric.

Now, $I_f$ contains two elements, corresponding to the orbit $\{ 1 \}$ and the orbit $\{\sigma_1, \sigma_1^{-1}\}$.\\

Another action of $G = \{1, -1\}$ on $A$ is given by $-1(\sigma) = -\sigma$,  ~i.e. when $G \subset \C^* \subset \aut_\C(A)$. This is an instance of the following example. \\

\noindent 2.  Let $A = \C[\sigma, \sigma^{-1}]$, and let $G = \{1, \zeta, \cdots, \zeta^{d-1} \}$, the group of $d$-th roots of unity, act on $A$ by $\zeta(\sigma^x) = \zeta^x \sigma^x$ ($G$ is a finite subgroup of $\C^* \subset \aut_\C(A)$). The one dimensional subspace $<\sigma^x>$ of $A$ spanned by $\sigma^x$ is $G$-invariant for each $x$ in $\mathbb{L}$. Let $\mathbb{L}_d = \{0, \pm d, \pm 2d \ldots \}$ be the sublattice of $\mathbb{L}$ defined by $d > 0$. If $d \nmid x$, then the trace of the map $\zeta^j: <\sigma^x> ~\rightarrow ~<\sigma^x>$ equals $\zeta^{xj}$, and hence the subspace $<\sigma^x>^G$ of fixed points of $<\sigma^x>$ has dimension $\frac{1}{d}(\chi(1) + \cdots + \chi(\zeta^{d-1})) = 0$. On the other hand, if $x = id$, then each map $\zeta^j: <\sigma^{x}> ~\rightarrow ~<\sigma^{x}>$ is the identity, and hence $\dim(<\sigma^{x}>^G) ~= ~1$. 
It follows that the dimension of $W_i^G$ equals $2i+1$, where $W_i$ is the subspace of $A$ spanned by $\{\sigma^x,~ -id \leq x \leq id \}$).

Let $I = (1 - \sigma^d)$, it is a $G$-invariant ideal of $A$. The dimension of $I_i^G = (I \cap W_i)^G$ equals $2i$ (it is spanned by $\{\sigma^{xd} - \sigma^{(x+1)d} | -r \leqslant x \leqslant r-1\}$), and thus $\dim((W_i/I_i)^G) = 1$. Exhausting $A$ by the nested sequence of finite dimensional $G$-invariant subspaces $W_i$, it follows that the dimension of $G$-invariant solutions $(A/I)'^G$ to $I$ equals 1. It is spanned by the solution $f$ which is 1 at points in $\mathbb{L}_d$ and 0 otherwise.

Indeed, by Proposition 4.2, the dimension of the space of $G$-invariant solutions equals $\frac{1}{d}\dim(A/I) = 1$. \\

\noindent 3. Let again $A = \C[\sigma, \sigma^{-1}]$, and now let $G = \{(1,1) (-1,-1)\} \subset \aut_\C(A)$ ($G$ is not a subgroup of either $\glo$ or $\C^*$ as in 1 or 2 above). $G$ acts on $A$ by 
$(-1,-1)(\sigma) = -\sigma^{-1}$. The ideal $I = (\sigma + \sigma^{-1})$ in 1(i) above is again $G$-invariant. The dimensions of $W_i$ and $I_i$ are unchanged; however now trace of $(-1,-1): I_i \rightarrow I_i$ equals -1, for all $i$. Thus, $\dim(I_i^G) = i-1$, and $\dim((W_i/I_i)^G) = 2$, for all $i$. Hence, the space $(A/I)'^G$ of $G$-invariant solutions has dimension 2. Now both $\cdots 0, -1, 0,\hat{1},0, -1, 0 \cdots$ and its shift $\cdots 1, 0, -1, \hat{0}, 1, 0, 1 \cdots$ are $G$-invariant solutions of $I$.      
\\

\noindent 4. Let $n = 2$, so that $\mathbb{L} = \Z^2$ and $A = \C[\sigma_1, \sigma_1^{-1}, \sigma_2, \sigma_2^{-1}]$. Let $G = S_2 = \{1, \tau\} \subset \glt$ be the group of permutations of $\{\sigma_1, \sigma_2\}$. This is the group of reflections of the lattice about the $\sigma_1 = \sigma_2$ diagonal.

Exhaust $A$ by the nested sequence of $G$-invariant finite dimensional subspaces $\{W_i\}_{i \geqslant 0}$, where $W_i$ is the subspace spanned by $\{\sigma_1^r \sigma_2^s ~|~ |r| + |s| \leqslant i \}$. A routine calculation shows that $\dim(W_i) = 2i^2 + 2i + 1$.
The map given by the transposition $\tau = (1,2)$ on $W_{2i}$, as well as on $W_{2i+1}$, has trace $2i + 1$; hence $\dim(W_{2i}^G) = 4i^2 + 3i + 1$ and $\dim(W_{2i+1}^G) = 4i^2 + 7i + 3$. Thus there is an infinite dimensional space of functions on $\Z^2$ that is symmetric with respect to $S_2$.

The ideal $I = (\sigma_1 - \sigma_2)$ is $G$-invariant. Let $I_i = I \cap W_i$. Its dimension equals $2i^2$; a basis is $\{\sigma_1^r \sigma_2^s - \sigma_1^{r-1}\sigma_2^{s+1} ~|~ |r| + |s| \leqslant i \}$. The trace of $\tau$ on $I_{2i}$ equals $-2i$, hence $\dim(I_{2i}^G) = \frac{1}{2}(\chi(1) + \chi(\tau)) = 4i^2 - i$. Similarly, the trace of $\tau$ on $I_{2i+1}$ equals $-2i - 2$, and $\dim(I_{2i+1}^G) = 4i^2 +3i$. Thus $\dim((W_{2i}/I_{2i})^G) = 4i + 1$ and $\dim((W_{2i+1}/I_{2i+1})^G) = 4i + 3$, or briefly $\dim((W_i/I_i)^G) = 2i + 1$ for all $i$. This is also the dimension of $W_i/I_i$, and thus every element in it is $G$-symmetric.

It follows that $(A/I)^G = A/I$, and hence that {\it every} solution in the infinite dimensional space of solutions to the equation defined by $(\sigma_1 - \sigma_2)$, is $G$-symmetric. Indeed the space of solutions coincides with the space of functions on $\mathbb{L}$ which are constant along the `anti-diagonal', and on each of its parallel translates. These functions are all $G$-symmetric.\\

Examples 1(ii), 3 and 4 point to the somewhat curious phenomenon where every solution to a $G$-invariant system $P$ of equations is $G$-invariant. Then, an arbitrary solution $f \in (A^k/P)'$ would need to satisfy $f([(q_1, \ldots, q_k]) = f([(gq_1, \ldots, gq_k)])$ for every $[(q_1, \ldots, q_k)] \in A^k/P$ and every $g \in G$. By Proposition 4.2 and Theorem 4.2, it is necessary that $G \subset \gl \subset \aut_\C(A)$, when its action restricts to the lattice $\mathbb{L}$.  As the points of $\mathbb{L}$ is a vector space basis for $A$, this translates to the following.

\begin{proposition} Let $G$ be a finite subgroup of $\gl$, and $P \subset A^k$ a $G$-invariant system of equations. Then every solution to $P$ is $G$-invariant if and only if $(\sigma^{x^1} - g\sigma^{x^1}, \ldots, \sigma^{x^k} - g\sigma^{x^k}) \in P$ for every $(\sigma^{x^1}, \ldots, \sigma^{x^k}) \in \mathbb{L}^k$ and every $g \in G$.
\end{proposition}	

In Example 1(ii), every $\sigma^x - \sigma^{-x}$ is in $J$, and in Example 3, every $\sigma_1^{x_1}\sigma_2^{x_2} - \sigma_1^{x_2}\sigma_2^{x_1}$ is in $I$. Thus both the ideals satisfy the condition of the above proposition.

\vspace{2mm}

Interesting questions on perturbations follow from 
the above results. We have shown that every system that is defined by a $G$-invariant submodule $P \subsetneq A^k$ admits a positive dimensional space $\ker_\F(P)^G$ of $G$-symmetric solutions. We could now construct a topology on the set $\mathcal{S}$ of all  systems, a (coarse) topology given as the direct limit of Zariski topologies on a filtration $\{\C^N\}$ of $\mathcal{S}$, as in \cite{sh}, with respect to which we would expect the dimension of $\ker_\F(P)^G$ to remain constant on algebraic subsets of $\S$. Such `stability' results are important in questions related to constructing models of  phenomena, because the coefficients of the monomials defining the system can only be approximately known.

\section{acknowledgement} I am indebted to Ananth Shankar for his generous help with this paper. I thank Debasattam Pal for many useful conversations, and the Electrical Engineering Department at IIT Bombay for its hospitality during many visits.


\begin{thebibliography}{9} 
	
	\bibitem{a} V.I.~Arnold, Mathematical Methods of Classical Mechanics, Graduate Texts in Mathematics 60, Springer-Verlag, 1978.
	
	\bibitem{d} I.V.~Dolgachev, Classical Algebraic Geometry: a modern view, Cambridge University Press, 2012.
	
	
	\bibitem{fw} F.~Fagnani and J.C.~Willems, Representations of symmetric linear dynamical systems, SIAM J. Control and Optimization, 31: 1267-1293, 1993.
	 
	\bibitem{mfg} M.~Folly-Gbetoula, Symmetry, reductions and exact solutions of the difference equation $u_{n+2} = au_n/(1+bu_nu_{n+1})$, Journal of Difference Equations and Applications, 23:1017-1024, 2017.
	 
	 \bibitem{hhs} M.~H\"{a}rk\"{o}nen, J.~Hirsch and B.~Sturmfels, Making waves, La Matematica, 2, 593-615, 2023. 
	 
	
	\bibitem{debs} D.~Pal and S.~Shankar, The coarsest lattice that determines a discrete multidimensonal system, Mathematics of Control, Signals, and Systems, 34: 405-433, 2022.    
	
	\bibitem{lp} L.~Peng, Relations between symmetries and conservation laws for difference systems, Journal of Difference Equations and Applications, 20:1609-1626, 2014.
	
	\bibitem{pom} J-F.~Pommaret, Partial differential control theory, volume II: control systems, Kluwer Academic Publishers, 2001.
	
	\bibitem{ser} J-P.~Serre, Linear representations of finite groups, Graduate Texts in Mathematics 42, Springer-Verlag, 1977.
	
	
	\bibitem{sh} S.~Shankar, The Hautus test and genericity results for controllable and uncontrollable behaviors, \phantom{x}SIAM jl. Control and Optimization, 52:32-51, 2014.
	
	\bibitem{stek} S. Shankar, Controllability and vector potential:~Six lectures at Steklov, \phantom{i} ~https://arxiv.org/abs/1911.01238, 2019.
	


\end{thebibliography}
\end{document}